\documentclass[noinfoline]{imsart}
\usepackage{jr,apajr,comment,enumerate,euscript,scalefnt,myfigures}
\newtheorem{metatheorem}{Meta-theorem}[section]

\arxiv{math.PR/0000512}

\begin{document}

\begin{frontmatter}

\title{On the trasductive arguments in statistics
}
\runtitle{On Transduction}

\begin{aug}
  \author{Ya'acov Ritov\thanksref{t2}%
  \ead[label=e1]{yaacov.ritov@gmail.com},
  \ead[label=u1,url]{http://pluto.mscc.huji.ac.il/$\sim$yaacov}}

  \thankstext{t2}{Research supported by an ISF grant.}

  \runauthor{Ya'acov Ritov}

  \affiliation{The Hebrew University of Jerusalem}

  \address{Department of Statistics, The Hebrew 
University, 91905 Jerusalem, Israel \printead{e1,u1}}

\end{aug}

\begin{abstract}
The paper argues that a part of the current statistical 
discussion is not based on the standard firm foundations 
of the field. Among the examples we consider are 
prediction into the future, semi-supervised 
classification, and causality inference based on 
observational data.
\end{abstract}

\begin{keyword}
\kwd{Foundations}
\kwd{Time series}
\kwd{Causality}
\kwd{Counterfactual}
\kwd{Semi-supervised learning}
\end{keyword}

\end{frontmatter}




\def\st{statistic\ }
\def\sts{statistics\ }
\def\stn{statistician\ }
\def\t{^{\scriptstyle{\mathsf{T}}} \hspace{-0.00em}}
\def\eq{}
\def\O{\EuScript{O}}
\def\sfa{\scalefont{0.6}}
\def\sfb{\scalefont{1.6666666666}}
\def\o{\ensuremath{\text{\sfa$\O$\sfb}}}

\section{introduction}
\label{sec:introduction}

Let $Y_1,\dots,Y_T$ be some time series. At time $T$ we 
want to predict the value of $Y_{T+1}$. This looks like a 
standard statistical problem,  feasible under an 
assumption of enough stationarity in the sequence. For 
example,  it may be assumed that 
$\eps_t=Y_t-\beta_1Y_{t-1}-\beta_2Y_{t-2}$, $t=3,4,\dots$ 
are \iid\  The value  $Y_{T+1}$ is going to be predicted 
by $\hat\beta_{T1}Y_T+\hat\beta_{T2}Y_{T-1}$,  where  
$\hat\beta_{T1}, \hat\beta_{T2}$  are estimated based on 
the sequence  $Y_1,\dots,Y_T$.  Does this practical 
thinking have a good statistical foundations?

Another example. Let $(X_1,Y_1), \dots, (X_n,Y_n) $ be an 
\iid sample, where $Y_i\in\{0,1\}$ is a label attached to 
observation $i$. Suppose we also have a large sample of 
unlabeled data $X_{n+1},\dots,X_N$, where $N\gg n$. Can 
we use the unlabeled data when we want to find a good 
classification rule? Can we justify this algorithm?

Finally, let $(X_1,Y_1), \dots, (X_n,Y_n) $ be a simple 
random sample, and $(X,Y)$ another copy, where, for 
simplicity, $X_i\in\{0,1\}$. We want to test whether $X$ 
is the cause of $Y$. Meaning, if we enforce $X$ to be 0, 
then the distribution of $Y$ will be different than if 
$X$ will be manipulated to be 1. Can this test be 
devised?

These three examples are typical to what we consider as a 
transduction inference, and we believe that most likely 
this type of inference go beyond the legitimate 
boundaries of standard statistical theory. In these 
examples, the statistician is extrapolating outside the 
observed model, to make prediction based on ungrounded 
belief.

A remark. The term statistics, as used in this paper, and 
its different derivatives like statistician, are not 
restricted to the research done by members of departments 
of statistics, who were also students in such a 
department. We mean by this anything related to inference 
done on the basis of empirically collected data, and by 
scientists that are dealing with anything from machine 
learning to neuroscience.

In the next section we will lay some of the standard 
foundation of the statistical practice. In Section 
\ref{sec:transduction} we explain what we mean by 
transduction. In the following three sections we describe 
in detail the three examples mentioned above. Some issues 
in time series are discussed in  Section 
\ref{sec:timeseries}. Inference with partially labeled 
data is described in Section \ref{sec:semisupervising}, 
and  Section \ref{sec:causality} deals with the causality 
argument. Some concluding remarks are given in Section 
\ref{sec:conclusions}.

\section{Background: the theoretical foundation of 
statistic inference}
\label{sec:background}

Statistical inference is based on an experiment. The 
basic notion of statistics is based on the collection of 
data, the understanding that at least in some sense the 
data could be different, and then, finally, making a 
statement which goes beyond the description of the actual 
observed data. Formally, the (statistical) experiment $E$ 
is built out a few elements. There is  a set $\Omega$, 
endued with a sigma field $\cal F$ and, a random element 
$Z$ which is measurable $\cal F$. There is a parameter 
set $\Theta$ and a family of probability measures 
$\{P_\theta: \;\theta\in\Theta\}$, such that for every 
$\theta\in\Theta$, $(\Omega,{\cal F},P_\theta)$ is a 
probability space.  We observe $Z$ and assume that it 
follows $P_\theta$  for some $\theta\in\Theta$.  See, for 
example, \cite{BW}.

The parameter $\theta$ is unknown, and is not directly 
observed. Some may call it the state of nature   
(\citeshortc{Berger1993}). In some very restrictive sense 
they are right.  It certainly is some parameter of 
reality. Whatever it is, we want to infer something about 
it from the observable $X$.  The statistician is 
reporting the evidence  about $\theta$ arising from the 
experiment $E$ and $Z$ (\citeshortc{Birnbaum1962}, 
\citeshortc{BW}).

Different authors differ on the scope of the statical 
conclusion. \cite{LeCam1986} believes that ``each element 
of the set \Th\ represents a particular `theory' about 
the physical phenomena involved in the experiment.'' In 
contrast, \cite{BickelDoksum} write ``Our aim is to use 
the data inductively, to narrow down in useful ways our 
ideas on what the `true' $P$ is.''  The difference 
between these two points of view may not seem apparent. 
We believe they answer differently the question of how 
far from the data the statistician can go. Le Cam speaks 
about a  theory underlined the data, and hence the 
conclusion from the experiment may go very far from the 
data, as far as the theory reaches, while Bickel and 
Doksum believe the statistician is limited to  
generalizations from sample to population. Thus, 
\cite{Cox1958} argues that ``a statistical inference 
carries us from observations to conclusions about the 
populations sampled.'' He contrasts this with ``a 
scientific inference in the broader sense [which] is 
usually concerned with arguing from descriptive facts 
about populations to some deeper understanding of the 
system under investigation.'' However, the leap from data 
to the deep understanding of Le Cam's theory about 
Berger's state of nature, is not statistical, or 
empirical, and most likely needs a leap of faith.

Thus, \cite{Tukey1960} points to the  ``difference 
between `statistical conclusions' and  `experimenter's 
conclusions'.''  Statisticians aim at precise statements, 
hence Tukey continues and claim that ``Both the 
statistician morale and integrity are tested \ldots when 
he has to face the possibility of a really substantial 
systematic error just after he has used all his skill to 
reduce, \ldots the effects of fluctuating errors to 95\% 
of their former value.''
The difficulty is when we go beyond the population, the 
safety net we create and proud of, like precise 
confidence interval, and P-values, are in doubt.

One of the basic concepts  of the field is that any 
precise statement on \th\ is impossible. Statistical 
inference is done with error. In other words, a 
\emph{particular inference} is rarely valid. Still, the 
field is proud in being able to make precise statements 
about the error. Whether this is done by presenting the 
object of inference as a random variable with a Bayesian 
\emph{a-posteriori} distribution, or with a frequentist 
confidence interval, the result is a precise 
quantification of the inference error. However, if the 
conclusion is derived with unknown  `systematic error', 
then one may doubt the importance of the exact 
quantification of the  'statistical error'.

Statisticians are well aware of the danger of 
extrapolation. For a very cute example consider the 
prediction of a newborn ear length. \cite{AltmanBland} 
used the regression line presented in 
\cite{Heathcote1995}: $\<\;\text{Ear\_Length} = 55.9 + 
0.22\times\text{Age}\;\>$ where ear length is measured in 
millimeter and age in years. This equation, based on a 
sample, predicts a ear length of 55.9mm at age 0, which 
is an absurd. The solution is simple: The regression line 
was found by measuring the ear length of a sample of 
30--93 years old men.  A minimum demands from a proper 
statistical inference is that it will be supported by the 
data. Extrapolation is going beyond the data, and hence 
it is considered problematic. In the next section we 
consider a much further reaching extrapolation, in which 
the statistician is going not only beyond the data, but 
also beyond the model.

\section{From induction and deduction to transduction}
\label{sec:transduction}

We argued that the legitimate statistical argument is 
from sample to population. For example, predicting the 
value of a random variable taken from the same 
distribution as the observed \iid\ sample. We refer to 
this type of statistical inference by induction. A 
statistical deduction is the inference about the 
parameter describing the population from which the sample 
was taken. The difference between these two may be 
considered verbal, pointing to two different perspectives 
on the same object. Consider for example a Gaussian 
mixture model: $X=(X_1,\dots,X_n)\in\R^n$, 
$P_\th=F_\th^n$, and $F_\th=\summ i1k \al_i 
\normal(\mu_i,1)$, where , $\th=(\al,\mu)$: $\al$ is a 
point in the $k$ dimensional simplex and $\mu\in\R^k$. 
The estimation of $F_\th(x)$ is what we call an 
induction, while the estimation of $\al$ is our 
deduction.

The problem of the justification of the statistical 
induction and deduction is on the one hand too simple to 
dwell on it, and on the other it is too deep and goes 
beyond the scope of this note.

Our concern is another type of statical inference, which 
we call transduction. When transduction is used, the 
statistician goes well beyond the data. This type of 
inference was relatively rare when the standard model was 
simple. If the statistician believe in the \iid shift 
normal model, there is very little that one can infer 
about except for the value of the mean (and maybe the 
variance) of the population distribution. When 
statistician started to work on more complex data 
structure, the possibilities and the dangers are 
abundant.

In fact, what is done here is extrapolating beyond the 
observed model, and not only beyond the range of the 
values observed in the sample. If the latter is 
dangerous, the former is much more so.

More formally. Let $\nu:\Th\to\scn$ be parameter of 
interest. Let $I_\nu(X;E)$ be any decision about 
$\nu(\th)$ done in the context of the experiment $E$. 
Statistician are trained to report the distribution of 
$I_\nu(X;E)$. Thus, the 95\% confidence interval is 
defined as a (random) set, which if the experiment will 
be repeated again and again, $X^i\dist F_{\th_i}$, 
$i=1,\dots,N$, and we will decide 
$I_\nu(X^1;E),\dots,I_\nu(X^N;E)$  then the cardinality 
of the set
$\{i: \th_i\in I_\nu(X^i;E)\}$ is approximately 0.95N. 
Taking care of the danger of extrapolation is a slightly 
more fuzzy. It means basically restricting the set of 
`legitimate' questions $\nu$ the statistician may ask, 
and at least avoiding  functions $\nu$ such that 
$\nu(\th)$ is not a smooth function of $F_\th$.

The \emph{gedankenexperiment} described above is what 
enables us to run away from the particular inference, on 
which we usually can say very little to the general 
scheme, which can be exact to a known degree.

However, in much of modern statistical theory, the 
\emph{gedankenexperiment} considered is different. In 
fact, the argument is based on  conceiving a list of 
experiment. $E^1,E^2,\dots$. The claim is, this type of 
logic worked in the past, why would it not work in this 
particular setup? Assuming that $P(Y|X)$ is smooth 
relative to the distribution of $X$, worked in this best 
typical (hardly related) models, why wouldn't it work for 
this new technique, for this very different model in a 
completely different field? We were able to prove that 
smoking is bad for your health, why wouldn't the proof 
that hormone replacement therapy is a wonderful cure of 
the pitfalls of the middle age be valid?

In the following sections we will discuss in some length 
the examples given in the introduction.

\section{Time series prediction}
\label{sec:timeseries}

Analyzing the role of the economists in the recent 
economic crises, Paul Krugman  wrote  ``the profession's 
blindness to the very possibility of catastrophic 
failures in a market economy. During the golden years, 
financial economists came to believe that markets were 
inherently stable — indeed, that stocks and other assets 
were always priced just right. There was nothing in the 
prevailing models suggesting the possibility of the kind 
of collapse that happened last year'' 
(\citeshortc{KrugmanNYT}). As he saw it, ``the economics 
profession went astray because economists, as a group, 
mistook beauty, clad in impressive-looking mathematics, 
for truth.'' We want to put this in a more general 
context (however, we do equate mathematical beauty with 
truth).

\subsection{Two problems of predictions}
We should distinguish between two very different 
problems.

\noindent\textbf{Prediction problem 1 (PP1):} Suppose 
$Y_1,\dots,Y_T,Y_{T+1}\in(-\en,\en)$. We observe 
$Y_1,\dots,Y_T$ and want to predict $Y_{T+1}$. A good 
predictor is such that
 \eqsplit{
    L(T)&=\bigl(\hat Y_{T+1}-  Y_{T+1}\bigr)^2
  }
is small.

\noindent\textbf{Prediction problem 2 (PP2):} Suppose 
that $Y_1,Y_2,\dots,Y_T \in [0,1]$. Suppose we want to 
predict $Y_{t+1}$ using its past, $Y_1,\dots,Y_t$, for  
$t=2,3\dots,T$. Suppose that a predictor is good if
 \eqsplit{
    L(T) &= \summ t1T (\hat Y_{t+1}-Y_{t+1})^2
  }
is small.

There are a few critical differences between these two 
problems. The first is that the first problem deals with 
a single event, while the second deals with a repeated 
one. As a result, the first problem asks for an 
unverified method, while the second problem asks for a 
method that can be checked and corrected as more data 
come in. The final difference between the two problems, 
is that the second has a built-in guarantee against 
catastrophe. The range of $Y_t$ was restricted to the 
unit interval. We believe that the first problem is less 
legitimate statistical problem than the second.

The main assumption underlined PP1 is that the future is 
in some sense  like the past. In some sense, since, 
typically we observe a dynamic system. This assumption, 
is not statistical.  Of course, PP1 is  well grounded if 
the prediction is based on a well verified physical 
theory. However, this is rarely the case. In the typical 
case, something like an ARMA model is going to be fitted 
to the existing data. Theory, if exists at all, is based 
on similar models used in the past, where they seemingly 
worked nicely. In either case we are in the 
$E^1,E^2,\dots$ \emph{Gedankenexperiment}, which prevents 
any possibility of giving sense to confidence intervals, 
or anything alike. The danger of the Black Swan phenomena 
is there, and it is real,  \cite{BlackSwan}.

PP2 is different. One is not interested only in the 
flock, and anyway, the range of their color is restricted 
from white to bright gray.  One needs very weak 
assumptions in order to solve it. In fact a very strong 
result can be claimed. We regress now into this model 
with a great detail. PP2 is directly related to the 
problem of self-calibrated deterministic forecaster, cf. 
\cite{Dawid1985}, \cite{Oakes1985}, \cite{Foster1999},   
\cite{FosterVohra}, and \cite{FudenbergLevine}.

\subsection{Regression. A solution to PP2}

One can start with any presumed estimator. E.g., the 
auto-regression model, and fit the parameters, then 
unbiased the prediction, given the past prediction 
experience. The result will be an estimator that predicts 
without significant bias the next observations, and 
preforms well along the sequence. This is so, without any 
assumption on the generator of $Y_1,Y_2,\dots$. In fact, 
the $Y$ sequence can be generated by somebody who knows 
the algorithm used by the statistician and tries to fool 
him. Here are some details.

Suppose first that we are observing the series 
$Z_1,Z_2,\dots$, and we want to construct an unbiased 
predictor such that $\E(Z_{t+1}  -  \hat Z_{t+1}\given 
\hat Z_{t+1}) =0$. To be more precise,  we consider a 
sort of zero-sum game in which the the statistician 
chooses functions $f_t:[0,1]^t\to[0,1]$, $t=1,2,\dots$. 
His opponent chooses a sequence $Z_1,Z_2,\dots$ (not 
necessarily randomly. The statistician predicts $Z_{t+1}$ 
by $\hat Z_{t+1}=f_t(Z_1,\dots,Z_t)$. Let $\eps>0$. The 
\stn wins if for any $z$, $|r_t(z) -z|>\eps $ a finite 
number of times where
 \eqsplit{
        r_t(z) &= \frac{\summ s1t Z_s\ind(|\hat 
Z_t-z|<\eps)}
        {\summ s1t Z_s\ind(|\hat Z_t-z|<\eps)} .
     }

The suggested predictor slowly walks on a grid according 
to a moving average of the observations. Let 
$\scs=\{\xi_0,\dots,\xi_K\}$, $\xi_j=2j\eta$, 
$\eta=1/2K$, be the
possible forecast values, and let $T$ be some large 
number, $T\gg
K$. Informally our procedure is as follows. If we switch
to the decision $\xi_j$ at the time $t_i$, we stick to 
this decision
for a while. We give the history that lead to $\xi_j$ a 
weight of
$T$, and accumulate new data. When the weighted mean 
deviates from
$\xi_j$ by more than $\eta$, we walk to  $\xi_{j\pm 1}$.

Formally, let $\phi_0\in\scs$ and $t_0=1$. For $i\geq 0$ 
define \eqsplit{
    t_{i+1} = \inf \Bigl\{t>t_{i}:\;
    \frac{T\phi_i+\summ s{t_i+1}t Z_s}{T+t-t_i}\not\in
    [\phi_i-\eta,\phi_i+\eta]\Bigr\},
}
where the infimum of an empty set is infinite, and
 \eqcases{\phi_{i+1}}{\phi_i+2\eta,
    \quad&
    \frac{T\phi_i+\summ s{t_i+1}{t_{i+1}} 
Z_s}{T+t_{i+1}-t_i} > \phi_i+\eta
    \\\\
    \phi_i-2\eta,&
    \frac{T\phi_i+\summ s{t_i+1}{t_{i+1}} 
Z_s}{T+t_{i+1}-t_i} < \phi_i-\eta,}
Note that $\phi_{i}\in\scs$, $i=1,2,\dots$. The suggested 
forecast
is $\hat Z_{t+1}=\phi_i$ for $t_i\leq t < t_{i+1}$.

Let \eqsplit{
    A_j&=\bigcup\bigl\{t_i,\dots,t_{i+1}-1 :\; 
\phi_i=\xi_j\;\&\; \phi_{i+1}=\xi_{j-1}\bigr\}
\\
    B_j&=\bigcup\bigl\{t_i,\dots,t_{i+1}-1: \; 
\phi_i=\xi_{j-1} \;\&\; \phi_{i+1}=\xi_j\bigr\},\quad 
j=1,\dots,
    K,
}
and
 \eqsplit{
    R_{tj} &= \frac{\summ s1t {Z_s}\ind(s\in A_j \union 
B_j)}
    {\summ s1t \ind(s\in A_j \union B_j)},\qquad 
t=1,2,\dots, j=1,\dots,K.
}
Then we have:

\begin{theorem}\label{th:r}
Under the above strategy, if all of  $t_1,t_2,\dots$ are 
finite,
then, for any $j=1,\dots,K$, either $\hat Z_{t+1}=\xi_j$ 
a finite number of
times, or $|R_{tj}-\xi _j|>2/K+K/T $ a finite number of 
times.
\end{theorem}

We return now to the the sequence $Y_1,Y_2,\dots$ of PP2. 
We start with $\ti Y_{t+1}$, a favorite predictor.  For 
example, in rain prediction $\ti Y_{t+1}$  can be based 
on any time series methods applied to the rain history, 
with or without the  information about the global weather 
map.  The next stage is approximately unbiased it. The 
canonical construction of the final predictor is based on 
a  partition $\{A_1,\dots,A_M\}$ of the possible values 
of $\ti Y_{t+1}$. For $m=1,\dots,M$ let 
$Z_{m1},Z_{m2},\dots$ be the subsequence $\{Y_t:\ti 
Y_t\in A_m\}$.  Apply the above scheme separately for 
each of these subsequences to obtain $\hat Z_{mt}$.  Then 
finally we use the predictor $\hat Y_{t+1}= \sum_m 
Z_{mT_m(t)}\ind(\ti Y_{t+1}\in A_m)$, where  
$T_m(t)=\#\{s:  s\le t \text{ and }\ti Y_{t+1}\in 
A_m\}$.

\subsection{Conclusion: the difference between PP1 and 
PP2}
The compound decision approach of PP2 saved the analysis, 
because it moved us from a particular  inference  (about 
$Y_T$, for q very particular $T$), to the  general 
inference  setup (about $Y_1,Y_2,\dots$). Also, we have 
restricted the possible values of the $Y$'s to avoid 
catastrophe by missing a single event.

We do not argue that prediction in time series is 
impossible. It is possible when it is a general scheme 
done under control (i.e., like PP2). We believe that 
predicting the future, that is, predicting one most 
important future event, is not a statistical task.

\section{Semi-supervising classification.}
\label{sec:semisupervising}

It is very frustrating to waste good data, even when  it 
is hardly related to the problem at hand. It is very 
tempting to use them, even if an unverified transductive 
argument is used in justifying the exercise.

\subsection{Can Classification be based on clustering?}

Consider a standard classification problem. A unit is 
characterized by a vector of variables  $X$ and a label 
$Y$. The statistical task is to find a classification 
rule by which the label $Y$ can be predicted based on the 
value of $X$.
It is quite standard that the data are collected 
electronically but judged by humans. It may then happen 
that most data points are not labeled.  This is the 
semi-supervised situation, and typical examples are 
satellite photography of terrains, or electronic 
espionage, but more down to earth examples exist, e.g., 
results of routine medical test may be in abundant, but a 
relatively few patients passed a through examination.

Formally, let $(X_1,Y_1), \dots, (X_n,Y_n) $ be an \iid 
sample, where $Y_i\in\{0,1\}$ is a label attached to 
observation $i$. Suppose we also have a large sample of 
unlabeled data $X_{n+1},\dots,X_N$, where $N\gg n$. In 
fact, $N$ may be large enough to be considered infinity 
for all practical purposes. See, for example, 
\cite{Joachims1999}, \cite{BelkinNiyogiRiemannian}, 
\cite{AminiGallinari}, and references given there.

\threefigures[1.3]{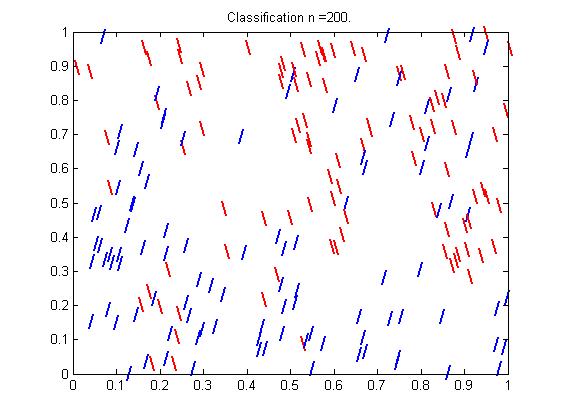}{transFig2.jpg}{transFig3.jpg}%
{(a) The classified data; (b) Unlabeled data was added; 
(c) The boundary of the Bayes classifier was 
added}{fig:semisupervised}

Semi-supervising classification is tempting. Consider 
Figure \ref{fig:semisupervised}(a). The problem is 
finding the best classifier. Geometrically, we want to 
find the boundary between the area were most of the 
slashes are normal to area where there are mostly the 
back slashes. It is not an easy task.  In Figure 
\ref{fig:semisupervised}(b) we added unlabeled points. 
Now, it seems not that difficult to find the boundary. 
Since the figures were created by simulation, we can add 
the optimal boundary corresponding to the Bayes 
classifier. It  is given in Figure 
\ref{fig:semisupervised}(c), lying  exactly where most 
readers thought it should be.

Here is an example where it is easy to see what is going 
on\footnote{This discussion was conceived in a discussion 
with Nicolai Meinshausen and Peter. B\"uhlmann }. Suppose 
$u_1,u_2,\dots,u_K$ are  \iid exponential with mean  
$\sig$, $v_1,v_2,\dots,v_K$ are  \iid exponential with 
mean $\tau$. Let $a_1<b_1<a_2<b_2<\dots<a_K<b_K$, be such 
that $b_i-a_i={u_i}/(\sum u_j+\sum v_j)$, and 
$a_{i+1}-b_i={v_i}/(\sum u_j+\sum v_j)$. Now let the 
covariate $X_1,X_2,\dots,X_N$ be \iid uniform on $\union 
(a_i,b_i)$, the variables $p_1,p_1,\dots,p_K$ be \iid 
$P(p_i=1)=P(p_i=0)=1/2$. Finally let the labels 
$Y_1,\dots,Y_n$ be independent $P(Y_i =1|X) = \sum p_j 
\ind\Bigl\{X\in (a_j,b_i)\Bigr\}$. Consider now the 
asymptotic as $n\to\en$, $n/N\to 0$, $n/K\to\gamma>0$,  
$\sig/\tau\to\en$, and $ N\tau/\log(N)K\to\en$.

The length of the support is 1, and there are $N$ 
observations uniformly distributed in it. Hence the 
largest spacing between observations belonging to same 
interval is $\O_p(\log(N)/N)$. The distance between 
adjacent intervals is approximately exponential, and 
hence the minimal distance between two adjacent intervals 
of the support is $\O_p(\tau /K )$. The ratio between 
these two is $\O_p\bigl(\log(N) K/N\tau\bigr)=\o_p(1)$. 
It follows, therefore, we can know where the spacing 
between any pair of adjacent unclassified observations is 
small and the two members  belong to the same interval, 
and where the spacing is large and they belong to 
different intervals.  Let $A_j$ and $B_j$ be the smallest 
and largest
observed $X$ inside the interval $(a_j,b_j)$. Then $\sum
|B_j-A_n|/\sum|b_j-a_j|\inc 1$. Hence with the 
unclassified data
it is easy to reconstruct the interval structure. If 
$\gamma$ is large, an almost perfect classifier can be 
constructed:  $\hat Y(x)=Y_{i(x)}$ where $i(x)$ is the 
first $i\in\{1,\dots,m_n\}$
such that $X_i$ and $x$ are in the same interval 
$[A_j,B_j]$, and
$1/2$ otherwise.

Let $(X_{(1)},Y_{(1)}),\dots,(X_{(n)},Y_{(n)})$ be the 
fully
observed sample sorted such that 
$X_{(1)}<X_{(2)}<\dots<X_{(n)}$.
Since $\sig/\tau\to \en$, most of the $(0,1)$ interval 
belongs to the support. Further,  $\gamma$ is finite, 
hence we cannot know whether $X_{(i)}$  and $X_{(i+1)}$ 
belong to the same interval or not. If $Y_{(i)}\ne
Y_{(i+1)}$, we know that $P(Y=1|x)$ is not constant on 
the interval $(X_{(i)}, X_{(i+1)}))$, but not much more 
than that. The change point is almost uniformly 
distributed (it is if $\gamma$ is very small). If  
$Y_{(i)}=
Y_{(i+1)}$, we know that it is quite likely that 
$P(Y=1|x)$ is not constant on the interval $(X_{(i)}, 
X_{(i+1)}))$. In any case we cannot know perfectly well 
$\E(Y|X)$ for a random $X$.

Clearly the classification error of an estimator based 
only on the unlabeled data is $1/2$---there is no way to 
know the values of the $p$'s given only the $X$'s. At the 
same time, any classifier based only on the label points 
is very weak and its classification error  is  close to 
the maximal  $1/2$. However, the classification error of 
an classifier based on all the data is close to the 
minimal 0. The semi-supervising approach works because 
the way $P(Y=1|x)$ is constructed is strongly tied to the 
way the support of $X$ is constructed, and the 
statistician knows these ties pretty well.

Can we justify the transduction from the experiment with 
observations $X$ to the experiment with $(X,Y)$? From the 
distribution of $X$ to the conditional distribution of 
$Y$ given $X$? Certainly the answer is yes, when 
simulations are done. However, can we give a statistical 
(empirical) justification for that? The answer to this 
difficulty is usually, something like ``see, it worked 
these many times, in this many best typical examples.'' 
We suggest to take these answers with the same grain of 
salt, as answers who dismiss the need of confidence sets, 
or P-values (or their Bayesian counterparts), because 
``rejecting the null was typically successful''.

The transduction argument was successful in the above two 
simulations, because, the data was generated by a 
mechanism that tied together the value of the regression 
function, and the underlined covariate distributions. 
Note, however the following pseudo-theorem. Cf. 
\cite{BKRW} for a precise formulations and examples.
\begin{metatheorem}
Suppose $(X,Y),(X_1,Y_1),(X_2,Y_2)\dots,$ are \iid. 
Suppose that  $X\dist H$, $P(Y=1|X=x)=p_{\th,\nu}(x)$, 
$(H,\nu)\in \sch\times\scn$. Then, the semiparametric 
bound for estimating $\th=\th(\nu)$ is the same whether 
the distribution $H$ of $X$ is known or unknown.
\end{metatheorem}
As we read the theorem it implies that there is no data 
dependent way to use the covariate distribution in a 
regular inference about the conditional distribution, at 
least locally and at the $\sqrt n$ rate. Any use of the 
covariate distribution, is based on \emph{a-priori} 
assumed connection and cannot be quantified or justified 
empirically.

\subsection{Preprocessing PCA}

Preprocessing PCA is another version of the 
semi-supervised transduction, cf. \cite{Jolliffe1972}, 
\cite{Jolliffe1973}, and \cite{CookForzani}. Suppose  
that $X,X_1,\dots,X_N$ are \iid, $X\in \R^p$. This sample 
is  partially labeled with a   $ Y\in\R$, such that 
$(X,Y),(X_1,Y_1,),\dots,(X_n,Y_n)$ are \iid, and $(X,Y)$ 
follows a linear regression model.  We consider the case 
where both $N\gg n$ and $p\gg n$. A common practice when 
there are too many variables, which is readily available 
in the semi-supervised state of mind, is reducing the 
number of variables using PCA (principal component 
analysis). This can be done using the all $X$ sample. 
After the PCA, we can regress $Y$ on the dominant $s$ 
main principal components, where $s\ll n$, and is chosen 
either $a-priori$, or in a data depended way (depending 
on $X_1,\dots,X_N$).

The logic is irresistible.  If $X=(X^1,\dots,X^p)\t$ and  
$\cor(X^1,X^2)\approx 1$, we certainly can retain only 
the average $(X^1+X^2)/2$, and ignore their difference. 
Thus reducing the number of variables by one. No much 
information is lost by this. However, this logic is 
transductive, and cannot be justified using the 
experimental $X$ data. It again based on some magic 
connection between the marginal distribution of $X$ and 
that of $Y$ given $X$. It appeals to some principal that 
says roughly that nothing is accidental. But this appeal 
to `justice' is clearly fallible. It may that   $Y$ is  a 
function of the mismatch between $X^1$ and $X^2$, and not 
so  much or their conjunction. For example, let $X^1$ be 
the optical power of the  lenses a patient needs (in 
diopters), and $X^2$ those he actually uses. Luckily they 
are highly correlated, but presumably headache is caused 
by their small difference. For another example, tension 
within a couple maybe caused more by education difference 
than by education average.

\subsection{Saving the transduction}

We have a real problem if by semi-supervised learning we 
understand what \cite{Gammerman98learningby} refer to by: 
``This is the problem of \emph{transduction} in the sense 
that we are interested in the classification of  a 
particular example rather than in the general rule for 
classifying future examples.''  However, with large data 
set a different point of view can be considered, and a 
more gentle interpretation is possible. Suppose we want 
to use the sample for prediction. The best predictor is 
nonparametric, and \emph{a-priori} belongs to a very 
large set of potential predictors, too large to make 
useless any empirical risk minimizer which is based on 
the unlabeled data. However, one may use the unlabeled 
sample essentially for suggesting a small set of 
potential predictors. The final predictor to be chosen, 
is going to be selected from these potential predictors, 
and this selection is going to be done solely based on 
the labeled sample. This predictor can be compared to the 
predictor which is based solely on the unlabeled data, 
and the best of them can be used. In this way, if the 
unverifiable assumptions used for the semi-supervising 
reasoning are approximately valid, then they will be 
utilized, and if they yield a bad predictor, they will be 
discarded.

Let us repeat. The unlabeled data  are used only in 
suggesting potential predictors, and not in the decision 
on the final predictor.

A similar argument can be used in the preprocessing PCA 
and the modified method is easily  defensible. 
Co-linearity is not really a problem for prediction. With 
LASSO like techniques (\citeshortc{Tibshirani1996}), 
$p\gg n$ can be handled, as long as there is an 
approximation with only $p_0\ll n$ non-zero coefficients, 
cf. \cite{GreenshteinRitov}. PCA preprocessing can be 
used to generate a new set of $q$, $q\ll p$ variables 
$Z^1,\dots, Z^q$. One can regress $Y$ on $Z^1,\dots,Z^q$, 
but for this one needs extra assumptions. One can regress 
$Y$ on $Z^1,\dots,Z^q,X^1,\dots,X^p$ with very little 
extra assumptions. Then if the PCA step happened to be 
smart, it will be effective, and if information was lost 
in the PCA reduction, it will be regained.

\section{Counterfactual causality }
\label{sec:causality}

The counterfactual theory of causality.(cf.  
\cite{Rubin1974}, \cite{Holland1988} and others)

\begin{itemize}

\item
Each individual is characterized by two outcomes 
$(Y^{C},Y^{T})$. One under the control condition and 
one under the treatment condition.

\item
The ``causal effect'' is the difference between these 
two potential outcomes, i.e., $\del=Y^{T} - Y^{C}$

\item
However,as mentioned, only one of these potential 
outcome is observed. The observation on subject $i$ 
is $(Y_i,T_i,X_i) $, where
 \eqsplit{
    Y_i=Y_i^C + T_i(Y_i^T-Y_i^C), \quad T\in\{0,1\}.
  }

\end{itemize}

For example, each participant carries two outcomes (from 
birth?), the first would be expressed if he will smoke 
all his life, and the other if he wouldn't. But the same 
subject may participate in another experiment, and 
therefore he has another couple of outcomes, where the 
first outcome will be measured  if he would learn German 
in one type of program, and the other will be expressed 
if he would learn it by another.

The model can be summarized by
$
    Y_i=Y_i^C + T_i\del_i, \quad T\in\{0,1\}
  $
When we are  dealing with a well designed experiment with 
a random allocations of units, $T$ is independent of 
$(Y^C,Y^T)$, and the mean causal effect is easily 
estimated
 \eqsplit{
    \hat\del &= \frac1n_1 \sum_{T_i=1}Y_i - \frac1n_0 
\sum_{T_i=0}Y_i
  }

However, this metaphysics is used exactly when $T$ is not 
exogenous and in particular it is  not  independent of 
$(Y^C,Y^T)$.
The different solutions of this basic ``difficulty'' are 
based on the assumption that $T$ is independent of 
$(Y^C,Y^T)$ conditional on a linear function of $X$.

Using a heavily loaded metaphysics in   a naturally 
positivistic science as statistics is justified when 
either
\begin{enumerate}
\item
It justifies in one sharp Ockhamian cut many many 
problems.
\item
It really simplifies the analysis.
\item
It unifies the terminology.
\end{enumerate}

Neither of these conditions is satisfied here. First, is 
the ax sharp?
Can the counterfactual theory of causality contribute 
anything a simple model cannot? What would be the case if 
the treatment is continuous? In reality, most 
``treatments'' are continuous even if measured as 
either-or. People are not either passive smoker or 
passive non-smoker, study to the exam or come unprepared, 
either take the drug or not. In medical experiment, even 
if the control condition is objectively defined, the 
experiment condition is typically arbitrary chosen from a 
continuous set of doses, treatment durations, Should we 
use a continuity of counterfactuals? What happens if the 
``treatment'' is multivariate? (Passive smoker in the 
work place, only once a hour, and once a week in 
pub\dots) A function of time? Simple it ain't.

Does the counterfactual point of view simplifies the 
terminology?
The natural terminology of the standard model $(\scz, 
\scf,\{P_\th\})$ is that of conditional distribution. 
Hence we ask, does the model $Y=Y^C+ T(Y^T-Y_C)$ has any 
additive value over: The conditional distribution of $Y$ 
given $T=0$ has a different mean than its conditional 
distribution given $T=1$?
We cannot dispense the conditional terminology altogether 
because we need to talk on the conditional distribution 
of $(T, Y^C,Y^T)$ given the covariate $X$. Hence, the 
counterfactual presentation does not simplifies the 
causality lexicon. It just adds new terms.

Some may say that it adds.
It tells us about two different statistical models: 
$(Y,\scf,\{Q^C_\th\}$ and $(Y,\scf,\{Q^T_\th\}$. These 
models are useful if we want to consider the distribution 
of $Y$ if the unit is ``enforced'' to be in one of the 
two groups. After all this is what we really want. We 
want to know what the impact of a new no-smoking rule 
will have on life expectancy.

However, if we consider  different statistical models, we 
could talk about a plentitude of them, not only on two. 
We can consider
$(Y,\scf,\{Q^t_\th\}$, $t\in\supp T$, and certainly it 
may that  $P_\th^t(Y\in A) = P_\th(Y\in A\given T=t)$, 
and for this we do not need the counter factual 
terminology.

This last presentation is a transduction over the 
standard statistical model which talks only about 
${P_\th}$. We don't have manipulated data in the sample, 
and hence any conclusion on manipulation is beyond the 
statistical reasoning. Typically it is just a wishful 
thinking. See for example \cite{BoundJJaegerBaker}, for a 
causality argument based self-evident assumptions that 
happened not to be true. Wishes sometimes come true, but 
more often than not, they do not. Thus \cite{Pearl2009} 
claims that ``confounding bias cannot be detected or 
corrected by statistical methods alone''. More 
specifically. ``This information must be provided by 
causal assumptions which identify relationships that 
remain invariant when external conditions change.'' But, 
most often than not, this is \emph{petitio principii}, or 
begging the question. One start with causal assumptions 
which are the empirical conclusion is disguise.

\section{Conclusions}
\label{sec:conclusions}

Transduction inference is unavoidable. One cannot avoid 
causality questions in the name of statistical integrity. 
One should realize that decision should be made, but one 
should recognize that they are not based on statistical 
safe ground. Certainly, quoting P-values and confidence 
intervals may be done only for the sake of style.

However, when we discussed the problems of prediction and 
semi-supervised learning we argued that when done 
carefully, almost transductive arguments can be used. In 
fact, one can use procedures that work nicely if his 
assumptions are valid, and do not fail him if they are 
not. If the time series is indeed stationary, then his 
predictors will be meaningful. If it is not, they would 
be trivial but right. The unlabeled data would be used to 
construct a better classifier if the right smoothness 
exists, otherwise it would not mislead the careful 
statistician.

\bibliographystyle{newapajr}  
\bibliography{transduction}

\end{document}